\documentclass[11pt,notitlepage,english]{article}
\usepackage{
  abstract,
  geometry,
  lmodern,
  mathtools,
  xparse
}

\usepackage[T1]{fontenc}
\usepackage[utf8]{inputenc}
\usepackage[sf,bf]{titlesec}
\usepackage[labelfont={bf,sf}]{caption}

\usepackage[amsmath,hyperref,thmmarks]{ntheorem}
\usepackage{amssymb,thmtools,upgreek}

\newcommand\proofSymbol{\ensuremath{_\blacksquare}}
\newcommand\qedSymbol{\proofSymbol}
\newcommand\qedhere{\ifmmode\qed\else\hfill\proofSymbol\fi}

\makeatletter
\def\theorem@checkbold{}
\newtheoremstyle{stmstyle}%
  {\item{\theorem@headerfont ##1\ ##2\theorem@separator}\enspace}%
  {\item{\theorem@headerfont ##1\ ##2\ {\normalfont(##3)}%
    \theorem@separator}\enspace}
\newtheoremstyle{prfstyle}%
  {\item{\theorem@headerfont ##1\theorem@separator}\enspace}%
  {\item{\theorem@headerfont ##3\theorem@separator}\enspace}
\makeatother
\theoremseparator{.}
\theorembodyfont{\upshape}

\theoremheaderfont{\itshape}
\theoremsymbol{\qedSymbol}
\declaretheorem[name=Proof,style=prfstyle]{prfenv}

\theoremheaderfont{\sffamily\bfseries}
\theoremsymbol{}

\declaretheorem[name=Proposition,style=stmstyle,numberlike=lemenv]%
  {proenv}

\usepackage[svgnames,x11names]{xcolor}
\usepackage[
  colorlinks=true,
  linkcolor=MediumBlue,
  citecolor=Firebrick2,
  hypertexnames=false
]{hyperref}
\usepackage[capitalise,compress,nameinlink]{cleveref}

\crefname{lemenv}{Lemma}{Lemmas}
\crefname{proenv}{Proposition}{Propositions}
\crefname{thmenv}{Theorem}{Theorems}
\crefname{corenv}{Corollary}{Corollaries}
\crefname{remenv}{Remark}{Remarks}
\crefname{exmenv}{Example}{Examples}
\crefname{figure}{Figure}{Figures}

\edef\envlabel{}
\newcommand\shortenv[1]{%
\DeclareDocumentEnvironment{#1}{ o g }{\ifvmode\else\unskip\fi%
  \IfValueTF{##1}{\begin{#1env}[##1]}{\begin{#1env}}%
  \IfValueT{##2}{\edef\envlabel{##2}\label{\envlabel}}}%
{\end{#1env}}%
}
\shortenv{lem}
\shortenv{pro}
\shortenv{thm}
\shortenv{cor}
\shortenv{rem}
\shortenv{exm}

\DeclareDocumentEnvironment{prf}{ o g }{%
  \IfValueTF{#1}{\begin{prfenv}[#1]}{\begin{prfenv}}}
{\end{prfenv}}

\newcommand\voc[1]{\emph{#1}}
\newcommand\fgn[1]{\emph{#1}}

\newcommand\bysame{\leavevmode\hbox to3em{\hrulefill}\thinspace}
\newcommand\MR[1]{\textsc{mr}:
  \href{http://www.ams.org/mathscinet-getitem?mr=#1}{\nolinkurl{#1}}}
\newcommand\ARXIV[1]{\textsc{arXiv}:
  \href{http://arXiv.org/abs/#1}{\nolinkurl{#1}}}

\usepackage{tikz}
\usetikzlibrary{calc}

\usepackage{mathrsfs}

\newcommand\NN{\mathbb N}
\newcommand\ZZ{\mathbb Z}
\newcommand\RR{\mathbb R}
\newcommand\CC{\mathbb C}
\newcommand\eps{\varepsilon}
\newcommand\es{\varnothing}
\usepackage{dsfont}
\newcommand\II{\mathds 1}

\newcommand\defeq{\coloneqq}
\newcommand\eqdef{\eqqcolon}
\newcommand\equiveq{\mathrel{\rlap{\raisebox{0.3ex}{$\cdot$}}%
  \raisebox{-0.3ex}{$\cdot$}}\equiv}
\newcommand\eqlaw{\stackrel d=}

\newcommand\ii{\mathrm i}
\newcommand\dd{\mathrm d}
\newcommand\bs[1]{\boldsymbol{\mathbf{#1}}}

\newcommand\limu{\mathop{\mathrm{lim}\!\!\uparrow}}
\newcommand\conv{\xrightarrow[n\to\infty]{}}
\newcommand\convlaw[1][n\to\infty]{\xrightarrow[#1]{(d)}}

\newcommand\oo[2]{\mathopen{(}#1,#2\mathclose{)}}

\newcommand\co[2]{\mathopen{[}#1,#2\mathclose{)}}

\usepackage{enumitem}
\newlist{assump}{enumerate}{1}
\setlist[assump]{label=\sffamily(H\arabic*)}
\crefname{assumpi}{Assumption}{Assumptions}

\newcommand\HT{\mathbb U}
\newcommand\HTh{\HT^{(h)}}

\newcommand\Ce{X}
\newcommand\Ceb{\bar X}
\newcommand\Cer{\hat X}
\newcommand\Cet{\widetilde X}
\newcommand\CeL{Y}
\newcommand\CeLb{\bar Y}

\newcommand\Hs{\mathtt A}
\newcommand\Sf[1][\le n\eps]{x^{#1}}
\newcommand\St[1][\le n\eps]{\tau^{#1}}
\newcommand\Stb[1][\le n\eps]{\bar\tau^{#1}}
\newcommand\StL{\bar\uptau^{\le\eps}}

\newcommand\GF{\bs X}
\newcommand\GFo{\mathds X}
\newcommand\GFL{\bs Y}

\newcommand\FT{\mathcal X}
\newcommand\FTo{\chi}
\newcommand\FTm{d_n}
\newcommand\FTL{\mathcal Y}

\DeclareMathOperator\height{ht}

\newcommand\Tk{p}

\newcommand\Ss{a}
\newcommand\Ssr[1]{A_{#1}}
\newcommand\Si{\gamma}

\newcommand\Jm{\Lambda}
\newcommand\Jmb{\bar\Lambda}
\newcommand\JmL{\Lambda}

\newcommand\Pa{{q_*}}
\newcommand\Pb{{q^*}}

\newcommand\Le{\Psi}
\newcommand\Cf{\kappa}
\newcommand\Cfb{\bar\Cf}
\newcommand\Cft{\widetilde\Cf}

\newcommand\Lp{\xi}
\newcommand\Lpb{\bar\xi}

\newcommand\Bt{\beta}
\newcommand\BtL{b}
\newcommand\Lt{\zeta}
\newcommand\Ltb{\bar\zeta}
\newcommand\LtL{\upzeta}
\newcommand\LtLb{\bar\upzeta}
\newcommand\Et{\mathcal E}
\newcommand\EtL{\epsilon}

\newcommand\Ba{M}

\renewcommand\Pr{\mathbb P}
\newcommand\Ex{\mathbb E}

\newcommand\Sk{\mathbb D}
\newcommand\leb[1][\Pb\downarrow]{\ell^{#1}}
\newcommand\Good{\mathcal G(n,\eps)}

\newcommand\RT{\mathscr T}
\newcommand\dRT{\mathrm{d}_{\mathrm{GH}}}

\title{\scshape Self-similar growth-fragmentations\\as scaling limits of
Markov branching processes}
\author{Benjamin Dadoun\thanks{Institut für Mathematik, 
  Universität Zürich, 
  Winterthurerstrasse 190, 
  CH-8057 Zürich, Switzerland.\hfill\eject
  Email: benjamin.dadoun@math.uzh.ch}}
\date{}
\begin{document}
\maketitle

\begin{abstract}
  We provide explicit conditions, in terms of the transition kernel of
  its driving particle, for a Markov branching process to admit a
  scaling limit toward a self-similar growth-fragmentation with negative
  index. We also derive a scaling limit for the genealogical embedding
  considered as a compact real tree.
\end{abstract}

\noindent{\itshape{\sffamily\bfseries Classification:} 60F17, 60J80}
\medskip
\section{Introduction}\label{sec:Newton}
Imagine a bin containing~$n$ balls which is repeatedly subject to random
(binary) divisions at discrete times, until every ball has been
isolated. There is a natural random (binary) tree with~$n$ leaves
associated with this partitioning process, where the subtrees above a
given height $k\ge0$ represent the different subcollections of all
$n$~balls at time~$k$, and the number of leaves of each subtree matches
the number of balls in the corresponding subcollection. The habitual
\voc{Markov branching property} stipulates that these subtrees must be
independent conditionally on their respective size. In the literature on
random trees, a central question is the approximation of so called
continuum random trees (CRT) as the size of the discrete trees tends to
infinity. We mention especially the works of
Aldous~\cite{Aldous91a,Aldous91b,Aldous93} and Haas, Miermont, et
al.\ \cite{Haas04,Miermont05,Haas08,Haas11,Haas12}.
Concerning the above example, Haas and Miermont~\cite{Haas12}
obtained, under some natural assumption on the splitting laws,
distributional scaling limits regarded in the
Gromov--Hausdorff--Prokhorov topology.
In the Gromov--Hausdorff sense where trees are considered as compact
metric spaces, they especially identified the so called
\voc{self-similar fragmentation trees} as the scaling limits. The latter
describe the genealogy of \voc{self-similar fragmentation processes},
which, reciprocally, are known to record the size of the components of
a (continuous) fragmentation tree above a given
height~\cite{Haas04}, and thus correspond to scaling limits for
partition sequences of balls as their number~$n$ tends to infinity. One
key tool in the work of Haas and Miermont~\cite{Haas12} is provided
by some non-increasing integer-valued Markov chain which,
roughly speaking, depicts the size of the subcollection containing a
randomly tagged ball. This Markov chain essentially captures the
dynamics of the whole fragmentation and, by their previous
work~\cite{Haas11}, itself possesses a scaling limit.

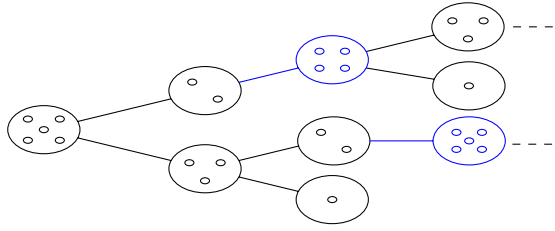
\begin{figure}[htp]
\centering
\begin{tikzpicture}[xscale=0.6,yscale=0.6]
  \draw[line width=.5pt] (0,0) circle (.8)
        (0,0) circle (.1)
        (45:.5) circle (.1)
        (135:.5) circle (.1)
        (225:.5) circle (.1)
        (315:.5) circle (.1)
        (20:.8) -- (20:3)
        (20:3.8) circle (.8)
        ($(30:.4)+(-20:3.8)$) circle (.1)
        ($(150:.4)+(-20:3.8)$) circle (.1)
        ($(-90:.4)+(-20:3.8)$) circle (.1)
        (-20:.8) -- (-20:3)
        (-20:3.8) circle (.8)
        ($(135:.4)+(20:3.8)$) circle (.1)
        ($(-45:.4)+(20:3.8)$) circle (.1)
        (20:7.6) -- (20:9.2)
        (20:10) circle (.8)
        ($(30:.4)+(20:10)$) circle (.1)
        ($(150:.4)+(20:10)$) circle (.1)
        ($(-90:.4)+(20:10)$) circle (.1)
        ($(20:6.8)+(-20:0.8)$) -- ($(20:6.8)+(-16:2.35)$)
        ($(20:6.8)+(-16:3.15)$) circle (.8)
        ($(20:6.8)+(-16:3.15)$) circle (.1)
        ($(-20:3.8)+(20:.8)$) -- ($(-20:3.8)+(20:2.2)$)
        ($(-20:3.8)+(18:3)$) circle (.8)
        ($(135:.4)+(-20:3.8)+(18:3)$) circle (.1)
        ($(-45:.4)+(-20:3.8)+(18:3)$) circle (.1)
        (-20:4.6) -- (-20:6)
        (-20:6.8) circle (.8)
        (-20:6.8) circle (.1);
  \draw[blue,line width=0.5pt]
        (20:4.6) -- (20:6)
        (20:6.8) circle (.8)
        ($(45:.4)+(20:6.8)$) circle (.1)
        ($(135:.4)+(20:6.8)$) circle (.1)
        ($(-45:.4)+(20:6.8)$) circle (.1)
        ($(-135:.4)+(20:6.8)$) circle (.1)
        ($(-20:3.8)+(18:3)+(0:.8)$) -- ($(-20:3.8)+(18:3)+(0:2.2)$)
        ($(-20:3.8)+(18:3)+(0:3)$) circle (.8)
        ($(45:.4)+(-20:3.8)+(18:3)+(0:3)$) circle (.1)
        ($(135:.4)+(-20:3.8)+(18:3)+(0:3)$) circle (.1)
        ($(-45:.4)+(-20:3.8)+(18:3)+(0:3)$) circle (.1)
        ($(-135:.4)+(-20:3.8)+(18:3)+(0:3)$) circle (.1)
        ($(-20:3.8)+(18:3)+(0:3)$) circle (.1)
        ($(-20:6.8)+(0:.8)$) -- ($(-20:6.8)+(0:2.2)$)
        ($(-20:6.8)+(0:3)$) circle (.8)
        ($(30:.4)+(-20:6.8)+(0:3)$) circle (.1)
        ($(150:.4)+(-20:6.8)+(0:3)$) circle (.1)
        ($(-90:.4)+(-20:6.8)+(0:3)$) circle (.1);
  \draw[dashed,line width=0.5pt] 
        ($(20:10)+(0:1)$) -- ($(20:10)+(0:2)$)
        ($(-20:3.8)+(20:3)+(-4:3)+(0:1)$) --
        ($(-20:3.8)+(20:3)+(-4:3)+(0:2)$)
        ($(20:6.8)+(-16:3.15)+(0:1)$) -- ($(20:6.8)+(-16:3.15)+(0:2)$)
        ($(-20:6.8)+(0:4)$) -- ($(-20:6.8)+(0:5)$);
\end{tikzpicture}
\caption{An example of dynamics with growth transitions (in blue).}%
\label{fig:Archimedes}
\end{figure}

The purpose of the present work is to study more general dynamics which
incorporate \emph{growth}, that is the addition of new balls in the
system (see \cref{fig:Archimedes}). One example of recent interest lies
in the exploration of random planar maps \cite{Kortchemski18,Curien18},
which exhibits ``holes'' (the yet unexplored areas) that split or grow
depending on whether the new edges being discovered belong to an already
known face or not. We thus consider a Markov branching system in
discrete time and space where at each step, every particle is replaced
by either one particle with a bigger size (growth) or by two smaller
particles in a conservative way (fragmentation). We condition the system
to start from a single particle with size~$n$ (we use the
superscript~$\cdot^{(n)}$ in this respect) and we are interested in its
behavior as $n\to\infty$. Namely, we are looking for:
\begin{enumerate}
\item A functional scaling limit for the process in
  time~$(\GFo(k)\colon k\ge0)$ of all particle sizes:
\begin{equation*}
  \left(\frac{\GFo^{(n)}(\lfloor\Ss_nt\rfloor)}{n}\colon t\ge0\right)
  \convlaw\,\bigl(\GFL(t)\colon t\ge0\bigr),
\end{equation*}
in some sequence space, where the $\Ss_n$ are positive (deterministic)
numbers;
\item A scaling limit for the system's genealogical tree, seen as a
  random metric space $(\FTo^{(n)},\FTm)$:
\begin{equation*}
  \left(\FTo^{(n)},\frac{\FTm}{\Ss_n}\right)
  \convlaw\,\FTL,
\end{equation*}
in the Gromov--Hausdorff topology.
\end{enumerate}
Like in the pure-fragmentation setting we may single out some specific
integer-valued Markov chain, but which of course is no longer
non-increasing. To derive a scaling limit for this chain, a first idea
is to apply, as a substitution to~\cite{Haas11}, the more general
criterion of Bertoin and Kortchemski~\cite{Kortchemski16} in terms of
the asymptotic behavior of its transition kernel at large states.
However, this criterion is clearly insufficient for the convergences
stated above as it provides no control on the ``microscopic'' particles.
To circumvent this issue, we choose to ``prune'' the system by freezing
the particles below a (large but fixed) threshold. That is to say, we
let the system evolve from a large size~$n$ but stop every individual as
soon as it is no longer bigger than some threshold $\Ba>0$ which will be
independent of~$n$, and we rather study the modifications~$\GF^{(n)}$
and~$\FT^{(n)}$ of the process and the genealogical tree that are
induced by this procedure.

The limits~$\GFL$ and~$\FTL$ that we aim at are, respectively, a
\voc{self-similar growth-fragmentation process} and its associated
genealogical structure. Indeed, the scaling limits of integer-valued
Markov chains investigated in~\cite{Kortchemski16}, which we build our
work upon, belong to the class of so called \voc{positive self-similar
Markov processes} (pssMp), and these processes constitute the
cornerstone of Bertoin's self-similar growth-fragmentations
\cite{Bertoin17,Curien18}. Besides, in the context of random planar maps
\cite{Kortchemski18,Curien18}, they have already been identified as
scaling limits for the sequences of perimeters of the separating cycles
that arise in the exploration of large Boltzmann triangulations.
Informally, a self-similar growth-fragmentation~$\GFL$ depicts a system
of particles which all evolve according to a given pssMp and whose each
negative jump~$-y<0$ begets a new independent particle with initial
size~$y$. In our setting, the self-similarity property has a negative
index and makes the small particles split at higher rates, in such a
way that the system becomes eventually extinct
\cite[Corollary~3]{Bertoin17}. The genealogical embedding~$\FTL$ is thus
a compact real tree; its formal construction is presented
in~\cite{Rembart18}.

Because of growth, one main difference with the conservative case is, of
course, that the mass of a particle at a given time no longer equals the
size (number of leaves) of the corresponding genealogical subtree. In a
similar vein, choosing the uniform distribution to mark a ball at random
will appear less relevant than a size-biased pick from an appropriate
(nondegenerate) supermartingale. This will highlight a Markov chain
admitting a self-similar scaling limit (thanks to the criterion%
~\cite{Kortchemski16}), and which we can plug into a many-to-one
formula. Under an assumption preventing an explosive production of
relatively small particles, we will then be able to establish our first
desired convergence. Concerning the convergence of the (rescaled)
trees~$\FT^{(n)}$, we shall employ a Foster--Lyapunov argument to obtain
an uniform control on their heights, which are nothing else than the
extinction times of the processes~$\GF^{(n)}$. Contrary to what one
would first expect, it turns out that a good enough Lyapunov function is
not simply a power of the size, but merely depends on the scaling
sequence~$(\Ss_n)$. This brings a tightness property that, together with
the convergence of ``finite-dimensional marginals'', will allow us to
conclude.

In the next section we set up the notation and the assumptions more
precisely, and state our main two results.

\section{Assumptions and results}\label{sec:Gauss}
Our basic data are probability transitions~$\Tk_{n,m},\,m\ge n/2$
and $n\in\NN$ ``sufficiently large'', with which we associate a Markov
chain, generically denoted~$\Ce$, that governs the law of the particle
system~$\GFo$: at each time $k\in\NN$ and with probability~$\Tk_{n,m}$,
every particle with size~$n$ either grows up to a size $m>n$, or
fragmentates into two independent particles with sizes 
$m\in\{\lceil n/2\rceil,\ldots,n-1\}$ and $n-m$. That is to say,
$\Ce^{(n)}(0)=n$ is the size of the initial particle in~$\GFo^{(n)}$,
and given~$\Ce(k)$ for some~$k\ge0$, $\Ce(k+1)$ is the largest among the
(one or two) particles replacing~$\Ce(k)$. We must emphasize that the
transitions~$\Tk_{n,m}$ from~$n$ ``small'' are irrelevant since our
assumptions shall only rest upon the asymptotic behavior of~$\Tk_{n,m}$
as~$n$ tends to infinity. Indeed, for the reason alluded in the
\hyperref[sec:Newton]{Introduction} that we explain further below,
we rather study the pruned version~$\GF$ where particles are frozen
(possibly at birth) when they become not bigger than a threshold
parameter~$\Ba>0$, which we will fix later on. Keeping the same
notation, this means that~$\Ce$ is a Markov chain stopped upon hitting
$\{1,2,\ldots,\Ba\}$. For convenience, \emph{we omit to write the
dependency in~$\Ba$}, and set $\Tk_{n,n}\defeq1$ for $n\le\Ba$.

In turn, the law of the genealogical tree~$\FT$ can be defined
inductively as follows (we give a more rigorous construction
in~\cref{sec:Euler}). Let $1\le k_1\le\cdots\le k_p$ enumerate the
instants during the lifetime~$\Lt^{(n)}$ of~$\Ce^{(n)}$ when
$n_i\defeq\Ce^{(n)}(k_i-1)-\Ce^{(n)}(k_i)>0$. Then~$\FT^{(n)}$ consists
in a branch with length~$\Lt^{(n)}$ to which are respectively attached,
at positions~$k_i$ from the root, independent trees distributed
like~$\FT^{(n_i)}$ (agreeing that~$\FT^{(n)}$ degenerates into a single
vertex for $n\le\Ba$). We view~$\FT^{(n)}$ as a metric space with metric
denoted by~$\FTm$.

Suppose $(\Ss_n)_{n\in\NN}$ is a sequence of positive real numbers which
is regularly varying with index~$\Si>0$, in the sense that for every
$x>0$,
\begin{equation}
  \lim_{n\to\infty}\frac{\Ss_{\lfloor nx\rfloor}}{\Ss_n}=x^\Si.%
  \label{eq:Riemann}
\end{equation}
Our starting requirement will be the convergence in
distribution
\begin{equation}
  \left(\frac{\Ce^{(n)}(\lfloor\Ss_nt\rfloor)}n\colon t\ge0\right)
  \convlaw\,\bigl(\CeL(t)\colon t\ge0\bigr),
\label{eq:Poincare}
\end{equation}
in the space $\Sk(\co{0}{\infty},\RR)$ of càdlàg functions
on~$\co{0}{\infty}$ (endowed with Skorokhod's J\textsubscript{1}
topology), towards a positive strong Markov process
$(\CeL(t)\colon t\ge0)$, continuously absorbed at~$0$ in an almost
surely finite time~$\LtL$, and with the following self-similarity
property:
\begin{flalign}
&\text{The law of~$\CeL$ started from $x>0$ is that of
$(x\CeL(x^{-\Si}\,t)\colon t\ge0)$ when~$\CeL$ starts from~$1$.}&%
\label{eq:Lagrange}
\end{flalign}
Since the seminal work of Lamperti~\cite{Lamperti72}, this simply means
that
\begin{equation}
  \log\CeL(t)=\Lp\!\left(\int_0^t\CeL(s)^{-\Si}\,\dd s\right)\!,\qquad
     t\ge0,%
  \label{eq:Euclid}
\end{equation}
with~$\Lp$ a Lévy process which drifts to~$-\infty$ as $t\to\infty$.
We denote by~$\Le$ the characteristic exponent of~$\Lp$ (so there is the
Lévy--Khinchine formula $E[\exp(q\Lp(t))]=\exp(t\Le(q))$ for every
$t\ge0$ and every~$q\in\CC$, wherever this makes sense), and by~$\JmL$
the Lévy measure of its jumps (that is a measure on~$\RR\setminus\{0\}$
with $\int(1\wedge y^2)\,\JmL(\dd y)<\infty$).

In order to state precisely our assumptions, we need to introduce some
more notation. First, we define the exponent
\begin{gather*}
  \Cf(q)\defeq\Le(q)+\int_{\oo{-\infty}{0}}\bigl(1-e^y\bigr)^q
    \,\JmL(\dd y),
\intertext{whose meaning will be discussed shortly (in the paragraph
``\hyperref[par:Moivre]{Discussion on the assumptions}''). Next, we also define,
for every~$n\in\NN$, the discrete versions}
  \Le_n(q)\defeq
    \Ss_n\sum_{m=1}^\infty\Tk_{n,m}
      \left[\left(\frac mn\right)^{\!q}-1\right]\!,\quad\text{and}\quad
  \Cf_n(q)\defeq\Le_n(q)+\Ss_n\sum_{m=1}^{n-1}\Tk_{n,m}
    \left(1-\frac mn\right)^{\!q}.
\end{gather*}
Finally, we fix some parameter~$\Pb>0$. After
\cite[Theorem~2]{Kortchemski16}, convergence~\eqref{eq:Poincare}
holds under the following two assumptions:
\begin{assump}
  \item\label{assump:H1} For every $t\in\RR$,
\begin{equation*}
  \lim_{n\to\infty}\Le_n(\ii t)=\Le(\ii t).
\end{equation*}
  \item\label{assump:H2} We have
\begin{equation*}
\limsup_{n\to\infty}\;
  \Ss_n\sum_{m=2n}^\infty\Tk_{n,m}\left(\frac mn\right)^{\!\Pb}
    <\,\infty.
\end{equation*}
\end{assump}
Indeed, by~\cite[Theorem~15.14 \&~15.17]{Kallenberg02}, \cref{assump:H1}
is essentially equivalent to (A1)\&(A2) of~\cite{Kortchemski16},
while~\labelcref{assump:H2} rephrases Assumption~(A3) there. We now
introduce the new assumption:
\begin{assump}[start=3]
  \item\label{assump:H3} We have
  either
  $\Cf(\Pb)<0$, or $\Cf(\Pb)=0$ and $\Cf'(\Pb)>0$. Moreover,
  for some $\eps>0$,
\begin{equation}
  \lim_{n\to\infty}
  \Ss_n\sum_{m=1}^{n-1}\Tk_{n,m}\left(1-\frac mn\right)^{\!\Pb-\eps}
  =\,\int_{\oo{-\infty}{0}}\bigl(1-e^y\bigr)^{\Pb-\eps}\,\JmL(\dd y).%
  \label{eq:H3}
\end{equation}
\end{assump}
Postponing the description of the limits, we can already state our two
convergence results formally:
\begin{thm}{thm:Hilbert}
  Suppose~\labelcref{assump:H1,assump:H2,assump:H3}. Then we can fix a
  freezing threshold~$\Ba$ sufficiently large so that, for every
  $q\ge1\vee\Pb$, the convergence in distribution
\begin{equation*}
  \left(\frac{\GF^{(n)}(\lfloor\Ss_nt\rfloor)}n\colon t\ge0\right)
  \convlaw\,\bigl(\GFL(t)\colon t\ge0\bigr),
\end{equation*}
  holds in the space $\Sk(\co{0}{\infty},\leb[q\downarrow])$,
  where~$\GFL$ is the self-similar growth-fragmentation driven
  by~$\CeL$, and
\begin{equation*}
  \leb[q\downarrow]\defeq
    \left\{\bs x\defeq(x_1\ge x_2\ge\cdots\ge0)\colon
      \sum_{i=1}^\infty(x_i)^q<\infty\right\}
\end{equation*}
(that is, the family of particles at a given time is always ranked in
the non-increasing order).
\end{thm}
\begin{thm}{thm:Grothendieck}
  Suppose \labelcref{assump:H1,assump:H2,assump:H3}, and $\Pb>\Si$. Then
  we can fix a freezing threshold~$\Ba$ sufficiently large so that there
  is the convergence in distribution
\begin{equation*}
  \left(\FT^{(n)},\frac{\FTm}{\Ss_n}\right)\convlaw\FTL,
\end{equation*}
  in the Gromov--Hausdorff topology, where~$\FTL$ is the random compact
  real tree that represents the genealogy of~$\GFL$.
\end{thm}

\paragraph{Description of the limits.}
As explained in the \hyperref[sec:Newton]{Introduction}, the
process~$\CeL$ portrays the size of particles in the self-similar
growth-fragmentation process~$\GFL\defeq(\GFL(t)\colon t\ge0)$, whose
construction we briefly recall (referring to~\cite{Bertoin17,Curien18}
for more details): The Eve particle~$\CeL_\es$ is distributed
like~$\CeL$. We rank the negative jumps of a particle~$\CeL_u$ in the
decreasing order of their absolute sizes (and chronologically in case of
\fgn{ex aequo}). When this particle makes its~$j$\textsuperscript{th}
negative jump, say with size $-y_j<0$, then a daughter
particle~$\CeL_{uj}$ is born at this time and evolves, independently of
its siblings, according to the law of~$\CeL$ started from~$y_j$. (Recall
that~$\CeL$ is eventually absorbed at~$0$, so we can indeed rank the
negative jumps in this way; for definiteness, we
set~$\BtL_{uj}\defeq\infty$ and~$\CeL_{uj}\equiveq0$ if~$\CeL_u$ makes
less than~$j$ negative jumps during its lifetime.) Particles are here
labeled on the Ulam--Harris tree~$\HT=\bigcup_{n=0}^\infty\NN^n$,
the set of finite words over~$\NN$, where $\NN^0=\{\es\}$ is reduced to the
root of the tree, and a vertex $u\defeq(u_1,u_2,\ldots,u_k)\in\HT$, at generation
$|u|\defeq k$, has $uj\defeq(u_1,u_2,\ldots,u_k,j)$ as $j$-th descendent.
Write~$\BtL_u$ for the birth time of~$\CeL_u$. Then
\begin{equation*}
  \GFL(t)=\Bigl(
    \CeL_u(t-b_u)\colon\,u\in\HT,\,\BtL_u\le t
  \Bigr),%
  \qquad t\ge0.
\end{equation*}
After~\cite{Bertoin17,Curien18}, this process is self-similar with
index~$-\Si$. Roughly speaking, this means that a particle with
size~$x>0$ evolves $x^{-\Si}$ times ``faster'' than a particle with
size~$1$. Since here~$-\Si<0$, there is the snowball effect
that particles get rapidly absorbed toward~$0$ as time passes, and it
has been shown \cite[Corollary~3]{Bertoin17} that such a
growth-fragmentation becomes eventually extinct, namely that
$\EtL\defeq\inf\{t\ge0\colon\GFL(t)=\emptyset\}$ is almost surely
finite.

The extinction time~$\EtL$ is also the height of the genealogical
structure~$\FTL$ seen as a compact real tree. Referring
to~\cite{Rembart18} for details, we shall just sketch the construction.
Let~$\FTL_{u,0}$ consists in a segment with length
$\LtL_u\defeq\inf\{t\ge0\colon\CeL_u(t)=0\}$
rooted at a vertex~$u$. Recursively, define $\FTL_{u,h+1}$ by attaching
to the segment~$\FTL_{u,0}$ the trees~$\FTL_{uj,h}$ at respective
distances $\BtL_{uj}-\BtL_u$ from~$u$, for each born particle~$uj$, $j\le h+1$.
The limiting tree $\FTL\defeq\limu_{h\to\infty}\FTL_h$, where
$\FTL_h\defeq\FTL_{\es,h}$ fulfills a so
called \voc{recursive distributional equation}. Namely, by
\cite[Corollary~4.2]{Rembart18}, given the sequence of negative jump
times and sizes $(\BtL_j,y_j)$ of~$\CeL$ and an independent sequence
$\FTL^1,\FTL^2,\ldots$ of copies of~$\FTL$, the action of grafting, on a
branch with length~$\LtL\defeq\inf\{t\ge0\colon\CeL(t)=0\}$ and at
distances~$\BtL_j$ from the root, the
trees~$\FTL^j$ rescaled by the multiplicative factor~$y_j^{\Si}$, yields
a tree with the same law as~$\FTL$. With this connection, Rembart and
Winkel \cite[Corollary~4.4]{Rembart18} proved that~$\EtL$ admits moments
up to the order $\sup\{q>0\colon\Cf(q)<0\}/\Si$. When particles do not
undergo sudden positive growth (i.e., $\JmL(\oo{0}{\infty})=0$), Bertoin
et al.\ \cite[Corollary~4.5]{Curien18} more precisely exhibited a
polynomial tail behavior of this order for the law of~$\EtL$.

\paragraph{Discussion on the assumptions.}\label{par:Moivre}
Observe that~\labelcref{assump:H1} entails~\labelcref{assump:H2} when
the Lévy measure~$\JmL$ of~$\Lp$ is bounded from above (in particular,
when~$\Lp$ has no positive jumps). By analyticity,
\cref{assump:H1,assump:H2} imply that $\Le_n(z)\to\Le(z)$ as
$n\to\infty$, for $0\le\Re z\le\Pb$. Adding the condition~\eqref{eq:H3}
in~\labelcref{assump:H3} then yields the convergence $\Cf_n(z)\to\Cf(z)$ for~%
$\Re z$ in a left-neighbourhood of~$\Pb$. Lastly, the first condition in~%
\labelcref{assump:H3} itself implies $\Le(\Pb)<0$ (since $\Le<\Cf$) and, together
with the other assumptions, that there must exist $\Pa\in\oo{0}{\Pb}$ and some
threshold~$\Ba$ such that
\[\Cf_n(q)<0\enspace\text{and}\enspace\Cf(q)<0,\qquad%
  \text{for all $q\in\co{\Pa}{\Pb}$\enspace and\enspace$n>\Ba$},\]
which is all but a superfluous requirement. Indeed, the condition $\Cf(q)\le0$ for
some $q>0$ is necessary (and sufficient)~\cite{Stephenson16} to prevent local
explosion of the growth-fragmentation~$\GFL$ (a phenomenon which would not allow us to
view it in some~$\ell^q$-space). Informally, the \voc{cumulant function}~$\kappa(q)$
captures the expected value of the sum of the particle sizes raised to the power~$q$
immediately after the first birth event. This function constitutes a key feature of
branching processes and, in particular, of self-similar growth-fragmentations~%
\cite{Shi17}. Of course, the meanings of the quantity~$\Cf_n(q)$ and of
the condition $\Cf_n(q)<0$ should be regarded the same but at the discrete level (that
is, w.r.t.~$\GF^{(n)}$).

\noindent We stress that our assumptions do
not provide any control on the ``small particles'' ($n\le\Ba$).
This explains why we need to ``freeze'' them (meaning that they
no longer grow or beget children); otherwise their number could
become quickly very high and make the system explode, as we illustrate
in the example below. We will basically choose~$\Ba$ as above, so that
$\Cf_n(q)\le0$ for some~$q$ and all~$n$, once we take the freezing into account (which
is tantamount to resetting%
\footnote{We make here a slight abuse on the notation. Again, even though the
dependency is not explicitly written, the discrete objects such as
$\Cf_n,\GF^{(n)},\ldots$ all ultimately depend on the freezing threshold~$\Ba$.}
$\kappa_n\equiveq0$ for $n\le M$).
\begin{exm}{exm:cex}
  Suppose that a particle with size~$n$ increases to size $n+1$ with
  probability $p<1/2$ and, at least when~$n$ is small, splits into
  two particles with sizes~$1$ and $n-1$ with probability $1-p$. Thus,
  at small sizes, the unstopped Markov chain essentially behaves like a
  simple random walk. On the one hand, we know from Cramér's theorem
  (see e.g.\ \cite[Theorem~2.2.3]{Dembo10})
  that for every~$\eps>0$ sufficiently small,
  \begin{equation*}
    \Pr\Bigl(\Ce^{(1)}(k)>(1-2p+\eps)k\Bigr),\qquad k\ge0,
  \end{equation*}
  decreases exponentially at a rate $c_p(\eps)>0$. On the other hand,
  keeping only track of particles with size~$1$ or~$2$,
  the number of particles with size~$1$ is bounded from below
  by~$Z^{[1]}$, where $\bs Z\defeq(Z^{[1]},Z^{[2]})$ is a $2$-type
  Galton--Watson process whose mean-matrix
  \begin{equation*}
    \begin{pmatrix}0&1\\2(1-p)&0\end{pmatrix}
  \end{equation*}
  has spectral radius $r_p\defeq\sqrt{2(1-p)}>1$, so that by the
  Kesten--Stigum theorem \cite[Theorem~V.6.1]{Athreya72}
  the number of particles with
  size~$1$ at time~$k\to\infty$ is of order at least~$r_p^k$, almost
  surely. Consequently the expected number of particles which are above
  $(1-2p+\eps)k$ at time~$2k$ is of exponential order at
  least $m_p(\eps)\defeq\log r_p-c_p(\eps)$. It is easily checked that
  this quantity may be positive (e.g., $m_{1/4}(1/4)\ge0.16$). Thus,
  without any ``local'' assumption on the reproduction law at small
  sizes, the number of small particles may grow exponentially and we
  cannot in general expect $\GFo^{(n)}(\lfloor\Ss_n\cdot\rfloor)/n$ to
  be tight in~$\leb[q\downarrow]$, for some $q>0$. However, this happens
  to be the case for the perimeters of the cycles in the branching
  peeling process of random Boltzmann
  triangulations~\cite{Kortchemski18}, where versions
  of~\cref{thm:Hilbert,thm:Grothendieck} hold for $\Si=1/2$, $\Pb=3$,
  and $\Ba=0$, although $\Cf_n(3)\le0$ seems fulfilled only for
  $\Ba\ge3$ (which should mean that the holes with perimeter~$1$ or~$2$
  do not contribute to a substantial part of the triangulation).
\end{exm}

We start with the relatively easy convergence of finite-dimensional
marginals (\cref{sec:Euler}). Then, we develop a few key results
(\cref{sec:Fermat21}) that will be helpful to complete the proofs of
\cref{thm:Hilbert} (\cref{sec:Hilbert}) and
\cref{thm:Grothendieck} (\cref{sec:Grothendieck}).

\section{Convergence of finite-dimensional marginals}\label{sec:Euler}
In this section, we prove finite-dimensional convergences for both the
particle process~$\GF$ and its genealogical structure~$\FT$. (We mention
that the freezing procedure is of no relevance here as it will be only
useful in the next section to establish tightness results; in particular
the freezing threshold~$\Ba$ will be fixed later.)

We start by adopting a representation of the particle system~$\GF$
that better matches that of~$\GFL$ given above. We define, for every
word $u\defeq u_1\cdots u_i\in\NN^i$, the
\voc{$u$-locally largest particle}%
\footnote{In the peeling of random Boltzmann maps~\cite{Kortchemski18},
the locally-largest cycles are called \voc{left-twigs}.}
$(\Ce_u(k)\colon k\ge0)$ by induction on $i=0,1,\ldots$ Initially,
for $i=0$, there is a single particle~$\Ce_\es$ labeled by $u=\es$,
born at time $\Bt_\es\defeq0$ and distributed like~$\Ce$. Then, we
enumerate the sequence $(\Bt_{u1},n_1),(\Bt_{u2},n_2),\ldots$ of the
negative jump times and sizes of~$\Ce_u$ so that $n_1\ge n_2\ge\cdots$
and $\Bt_{uj}<\Bt_{u(j+1)}$ whenever $n_j=n_{j+1}$. Conditionally on
$(n_j)_{n\ge1}$, the processes
$\Ce_{uj},\,j=1,2,\ldots$ are independent and distributed like
$\Ce^{(n_j)}$ respectively (for definiteness, we set
$\Bt_{uj}\defeq\infty$ and $\Ce_{uj}\equiveq0$ if~$\Ce_u$ makes less
than~$j$ negative jumps during its lifetime), and we have
\begin{equation*}
  \bigl(\GF(k)\bigr)_{k\ge0}
    \eqlaw\Bigl(\Ce_u(k-\Bt_u)\colon u\in\HT,\,
    \Bt_u\le k\Bigr)_{k\ge0}.
\end{equation*}
Recall the notation~$\cdot^{(n)}$ to stress that the system is started
from a particle with size $\Ce_\es(0)=n$.
\begin{lem}{lem:Euler}
  Suppose \labelcref{assump:H1,assump:H2}. Then for every finite subset
  $U\subset\HT$, there is the convergence
  in~$\Sk(\co{0}{\infty},\RR^U)$:
\begin{equation}
  \left(\frac{\Ce_u^{(n)}(\lfloor\Ss_nt\rfloor-\Bt^{(n)}_u)}n
    \colon u\in U\right)_{t\ge0}
    \convlaw\,\Bigl(\CeL_u(t-\BtL_u)\colon u\in U\Bigr)_{t\ge0}.%
  \label{eq:Poincare2}
\end{equation}
\end{lem}
\begin{prf}{prf:Euler}
  We follow the argument used to prove the second part of
  \cite[Lemma~17]{Kortchemski18}. For $h\ge0$, let
  $\HT_h\defeq\{u\in\HT\colon|u|\le h\}$ be the set of vertices
  with height at most~$h$ in the tree~$\HT$. It suffices to show
\begin{flalign*}
  &(\mathscr I_h)\colon
  &\text{Convergence~\eqref{eq:Poincare2} holds in
    $\Sk(\co{0}{\infty},\RR^U)$ for every finite subset
    $U\subset\HT_h$},&&&&&&&&&&
\end{flalign*}
  by induction on~$h$. The statement~$(\mathscr I_0)$ is given
  by~\eqref{eq:Poincare}. Now, if~$U$ is a finite subset
  of~$\HT_{h+1}$ and $F_u,\,u\in U,$ are continuous bounded functions
  from $\Sk(\co{0}{\infty},\RR)$ to~$\RR$, then the branching property
  entails that, for
  $\Cer_u^{(n)}\defeq\Ce_u^{(n)}(\lfloor\Ss_n\cdot\rfloor-\Bt_u)/n$,
\begin{equation*}
  \Ex\!\left[\prod_{u\in U}\!
    F^{\vphantom{(n)}}_u\bigl(\Cer^{(n)}_u\bigr)\;\bigg|\;
    \bigl(\Ce_u^{(n)}\colon u\in\HT_h\bigr)\right]
  \,\,=\!\!\prod_{u\in U\cap\,\HT_h}\!\!\!\!\!
    F^{\vphantom{(n)}}_u\bigl(\Cer^{(n)}_u\bigr)
    \;\cdot\!\!
    \prod_{\substack{u\in U\\|u|=h+1}}\!\!\!\!
      E_{\Cer^{(n)}_u(0)}^{(n)}\bigl[F_u\bigr],
\end{equation*}
  where~$E_x^{(n)}$ stands for expectation under the law~$P_x^{(n)}$
  of $\Cer_\es$ started from~$x$, which
  by~\eqref{eq:Riemann}, \eqref{eq:Poincare} and~\eqref{eq:Lagrange},
  converges weakly as $n\to\infty$ to the law~$P_x$ of~$\CeL$
  started from~$x$. The values $\Cer^{(n)}_u(0)$ for $|u|=h+1$
  correspond to (rescaled) negative jump sizes of particles at
  height~$h$. With~\cite[Corollary~VI.2.8]{Jacod87} and our convention
  of ranking the jump sizes in the non-increasing order, the
  convergence in distribution
  $(\Cer^{(n)}_u(0)\colon u\in U,\,|u|=h+1)
    \to(\CeL_u(0)\colon u\in U,\,|u|=h+1)$ as $n\to\infty$
  thus holds jointly with~$(\mathscr I_h)$. Further, thanks to the
  Feller property~\cite[Lemma~2.1]{Lamperti72} of~$\CeL$, its
  distribution is weakly continuous in its starting point. By the
  continuous mapping theorem we therefore obtain, applying back the
  branching property, that
\begin{equation*}
  \Ex\!\left[\prod_{u\in U}\!F^{\vphantom{(n)}}_u
    \bigl(\Cer^{(n)}_u\bigr)\right]
  \conv\,\Ex\!\left[\prod_{u\in U}\!
    F_u\bigl(\CeL_u(\cdot-\BtL_u)\bigr)\right]\!.
\end{equation*}
  A priori, this establishes the convergence in distribution
  $(\Cer^{(n)}_u\colon u\in U)
    \to(\CeL_u(\cdot-\BtL_u)\colon u\in U)$ only in
  the product space $\Sk(\co{0}{\infty},\RR)^U$. By
  \cite[Proposition~2.2]{Jacod87} it will also hold in
  $\Sk(\co{0}{\infty},\RR^U)$ provided that the processes
  $\CeL_u,\,u\in U,$ almost surely never jump simultaneously. But this
  is plain since particles evolve independently and the jumps of~$\CeL$
  are totally inaccessible. Thus
  $(\mathscr I_h)\implies(\mathscr I_{h+1})$.
\end{prf}

Next, we proceed to the convergence of the finite-dimensional marginals
of~$\FT$, which we shall first formally construct. For each $u\in\HT$
with $\Bt_u<\infty$, let~$\Lt_u$ denote the lifetime of the stopped
Markov chain~$\Ce_u$. Recall the definition in~\cref{sec:Gauss} of the
trees $\FTL,\FTL_h,\,h\ge0$, related to~$\GFL$'s genealogy, that echoes
Rembart and Winkel's construction~\cite{Rembart18}. Similarly,
let~$\FT_{u,0}$ simply consist of an edge with length~$\Lt_u$, rooted at
a vertex~$u$. Recursively, define~$\FT_{u,h+1}$ by attaching to the
edge~$\FT_{u,0}$ the trees~$\FT_{uj,h}$ at a distance $\Bt_{uj}-\Bt_u$
from the root~$u$, respectively, for each born particle~$uj,\,j\le h+1,$
descending from~$u$. The tree $\FT_h\defeq\FT_{\es,h}$ is a finite tree
whose vertices are labeled by the set~$\HTh$ of words over
$\{1,\ldots,h\}$ with length at most~$h$. Plainly, the sequence
$\FT_h,\,h\ge0,$ is consistent, in that~$\FT_h$ is the subtree
of~$\FT_{h+1}$ with vertex set~$\HTh$, and we may consider the inductive
limit $\FT\defeq\limu_{h\to\infty}\FT_h$. We write $\FTm(v,v')$ for the
length of the unique path between~$v$ and~$v'$ in~$\FT^{(n)}$. All these
trees belong to the space~$\RT$ of (equivalence classes of) compact,
rooted, real trees and can be embedded as subspaces of a large metric
space (such as, for instance, the space $\ell^1(\NN)$ of summable
sequences \cite[Section~2.2]{Aldous93}).
Irrespectively of the embedding, they can be compared one with each
other through the so called Gromov--Hausdorff metric~$\dRT$ on~$\RT$. We
forward the reader to \cite{LeGall06,Evans08} and references therein.
\begin{lem}{lem:Galois}
  Suppose \labelcref{assump:H1,assump:H2,assump:H3}. Then for all
  $h\in\NN$, there is the convergence in $(\RT,\dRT)$:
  \begin{equation*}
    \left(\FT^{(n)}_h,\frac{\FTm}{\Ss_n}\right)\convlaw\,\FTL_h.
  \end{equation*}
\end{lem}
\begin{prf}
It suffices to show the joint convergence of all branches. The branch
going from the root~$\es$ through the vertex~$u\in\HTh$ has total length
$\Et^{(n)}_u\defeq\Bt^{(n)}_u+\Lt^{(n)}_u$ in~$(\FT^{(n)}_h,\FTm)$, and
length $\EtL_u\defeq\BtL_u+\LtL_u$ in~$\FTL_h$. Recall that
conditionally on $\{\Ce_u(0)=n\}$, the random variable~$\Lt_u$ has the
same distribution as $\Lt^{(n)}\defeq\Lt^{(n)}_\es$.
By~\cite[Theorem~3.(i)]{Kortchemski16}, the convergence
\begin{equation*}
  \frac{\Lt^{(n)}}{\Ss_n}\convlaw\LtL\defeq\inf\{t\ge0\colon\CeL(t)=0\}
\end{equation*}
holds jointly with~\eqref{eq:Poincare}. Adapting the proof
of~\cref{lem:Euler}, we can more generally check that for every
finite subset~$U\subset\HT$, we have, jointly
with~\eqref{eq:Poincare2},
\begin{equation*}
  \left(\frac{\Et^{(n)}_u}{\Ss_n}\colon u\in U\right)
  \convlaw\bigl(\EtL_u\colon u\in U\bigr).
\end{equation*}
In particular, this is true for $U\defeq\HTh$.
\end{prf}

To conclude this section, we restate an observation of Bertoin, Curien,
and Kortchemski \cite[Lemma~21]{Kortchemski18} which results from the
convergence of finite-dimensional marginals (\cref{lem:Euler}): with
high probability as $h\to\infty$, ``non-negligible'' particles have
their labels in~$\HTh$. Specifically, say that an individual $u\in\HT$
is $(n,\eps)$-good, and write $u\in\Good$, if the particles~$\Ce_v$
labeled by each ancestor~$v$ of~$u$ (including~$u$ itself) have size at
birth at least~$n\eps$. Then:
\begin{lem}{lem:Neumann}
We have
\begin{equation*}
  \mathop{\vphantom{\limsup}\lim}_{h\to\infty}
  \limsup_{\vphantom{h}n\to\infty}\;
    \Pr^{(n)}\!\left(\Good\not\subseteq\HTh\right)=\,0.
\end{equation*}
\end{lem}

\section{A size-biased particle and a many-to-one formula}%
\label{sec:Fermat21}
We now introduce a ``size-biased particle'' and relate it to a
many-to-one formula. This will help us derive tightness estimates in
\cref{sec:Hilbert,sec:Grothendieck}, and thus complement the
finite-dimensional convergence results of the preceding section.
Recall from \cref{assump:H1,assump:H2,assump:H3} that we can find
$\Pa\in\oo{0}{\Pb}$ such that, as $n\to\infty$, $\Cf_n(q)\to\Cf(q)<0$ for
every $q\in\co{\Pa}{\Pb}$.
Consequently, \emph{we may and will suppose for the remainder of this
section that the freezing threshold~$\Ba$ is taken sufficiently large so
that $\Cf_n(\Pa)\le0$ for every $n>\Ba$}. (Note that $\Cf_n(\Pa)=0$
for $n\le\Ba$, by our convention $\Tk_{n,n}\defeq1$.)
\begin{lem}{lem:Abel}
  For every $n\in\NN$,
\begin{equation}
  \Ex^{(n)}\!\left[\bigl(\Ce(1)\bigr)^\Pa
    +\bigl(n-\Ce(1)\bigr)_+^\Pa
  \right]\le\,n^\Pa.%
  \label{eq:Abel}
\end{equation}
  Therefore, the process%
\footnote{We set here $\Ce_u(i)\defeq0$ for $i<0$ in order to not
burden the notation with the indicator $\II_{\{\Bt_u\le k\}}$.}
\begin{equation*}
\sum_{u\in\HT}
   \bigl(\Ce_u(k-\Bt_u)\bigr)^\Pa,\qquad k\ge0,
\end{equation*}
is a supermartingale under~$\Pr^{(n)}$.
\end{lem}
\begin{prf}{prf:Abel}
  The left-hand side of~\eqref{eq:Abel} is
\begin{equation*}
    n^\Pa\sum_{m=0}^\infty\Tk_{n,m}
      \left[\left(\frac mn\right)^{\!\Pa}
        +\left(1-\frac mn\right)_{\!+}^{\!\Pa}\right]
    =\,n^\Pa\left(1+\frac{\Cf_n(\Pa)}{\Ss_n}\right)\!,
\end{equation*}
  where $\Cf_n(\Pa)\le0$. Hence the first part of the statement. The
  second part follows by applying the branching property at any given
  time~$k\ge0$:
\begin{flalign*}
\Ex^{(n)}\!\left[\sum_{u\in\HT}\bigl(\Ce_u(k+1-\Bt_u)\bigr)^\Pa\;
   \Big|\;\GF(k)=(x_i\colon i\in I)\right]
&=\,\sum_{i\in I}\Ex^{(x_i)}\!\left[\bigl(\Ce(1)\bigr)^\Pa
    +\bigl(x_i-\Ce(1)\bigr)_+^\Pa\right]\\[.4em]
&\le\,\sum_{i\in I}(x_i)^\Pa\\[.4em]
&=\,\sum_{u\in\HT}
   \bigl(\Ce_u(k-\Bt_u)\bigr)^\Pa.
\end{flalign*}
\end{prf}
\begin{rem}
  Put differently, the condition ``$\Cf_n(\Pa)\le0$'' entails that
  $n\mapsto n^\Pa$ is superharmonic with respect to the
  ``fragmentation operator''. This map plays the same role as the
  function~$f$ in~\cite{Kortchemski18}, where it takes the form of a
  cubic polynomial ($\Pa=3$) and
  \[\sum_{u\in\HT}f\bigl(\Ce_u(k-\Bt_u)\bigr),\qquad k\ge0,\]
  is actually a martingale. More generally, the map
  $n\mapsto n^\Pa$ could be replaced by any regularly-varying sequence
  with index~$\Pa$, but probably at the cost of heavier notation.
\end{rem}
As we see in the proof of \cref{lem:Abel}, the fact that $\Cf_n(\Pa)\le0$ allows us to
introduce a (defective) Markov chain $(\Ceb(k)\colon k\ge0)$ on~$\NN$, to which we
add~$0$ as cemetery state, with transition
\begin{equation}
  \Ex^{(n)}\!\left[f\bigl(\Ceb(1)\bigr);\,\Ceb(1)\neq0\right]
    =\sum_{m=1}^\infty\Tk_{n,m}\left[\left(\frac mn\right)^{\!\Pa}
        f(m)
      +\left(1-\frac mn\right)_{\!+}^{\!\Pa}
        f(n-m)\right]\!.%
  \label{eq:Weierstrass}
\end{equation}
We let~$\Ltb\defeq\inf\{k\ge0\colon\Ceb(k)=0\}$ denote its lifetime.
Up to a change of probability measure, $\Ceb$ follows the trajectory of
a randomly selected particle in~$\GF$, until it is eventually absorbed to the
cemetery state~$0$. It admits the following scaling
limit (which could also be seen as a randomly selected particle
in~$\GFL$; see~\cite[Section~4]{Curien18}):
\begin{pro}{pro:Descartes}
  There is the convergence in distribution
\begin{equation}
  \left(\frac{\Ceb^{(n)}(\lfloor\Ss_nt\rfloor)}n\colon t\ge0\right)
  \convlaw\,\left(\CeLb(t)\colon t\ge0\right)
  \label{eq:Descartes}
\end{equation}
  in $\Sk(\co{0}{\infty},\RR)$, where the limit~$\CeLb$ fulfills the
  same identity~\eqref{eq:Euclid} as~$\CeL$, but for a (killed) Lévy
  process~$\Lpb$ with characteristic exponent~$\Cfb(q)\defeq\Cf(\Pa+q)$.
  Further, if~$\LtLb$ denotes the lifetime of~$\CeLb$, then the
  convergence
\begin{equation*}
  \frac{\Ltb^{(n)}}{\Ss_n}\convlaw\LtLb
\end{equation*}
  holds jointly with~\eqref{eq:Descartes}.
\end{pro}
\begin{prf}{prf:Descartes}
  Write~$\Jmb_n$ for the law of~$\log(\Ceb^{(n)}(1)/n)$, with the
  convention $\log0\defeq-\infty$. We see from~\eqref{eq:Weierstrass}
  that $\Ss_n\,\Pr(\Ceb(1)=0)=-\Cf_n(\Pa)$, and, for every
  $0\le q\le\Pb-\Pa$,
\begin{align*}
  \int_\RR
    (e^{qy}-1)\,\Ss_n\Jmb_n(\dd y)
  &=\Ss_n\sum_{m=0}^\infty\Tk_{n,m}\left[\left(\frac mn\right)^{\!\Pa}
      \left(\!\left(\frac mn\right)^{\!q}-1\right)
        +\left(1-\frac mn\right)_{\!+}^{\!\Pa}
    \left(\!\left(1-\frac mn\right)^{\!q}-1\right)\right]\\
  &=\Cf_n(\Pa+q)-\Cf_n(\Pa).
\end{align*}
Hence
\begin{equation*}
  -\Ss_n\Jmb_n(\{-\infty\})+\int_\RR
    (e^{qy}-1)\,\Ss_n\Jmb_n(\dd y)
  =\Cf_n(\Pa+q)\conv\Cfb(q).
\end{equation*}
Furthermore, by~\labelcref{assump:H2},
\begin{equation*}
  \limsup_{n\to\infty}\;\Ss_n\int_1^\infty e^{(\Pb-\Pa)y}
    \,\Jmb_n(\dd y)
  \le\limsup_{n\to\infty}\;\Ss_n\sum_{m=2n}^\infty\Tk_{n,m}
    \left(\frac mn\right)^{\!\Pb}<\,\infty.
\end{equation*}
In other words, assumptions~(A1), (A2) and~(A3) of~\cite{Kortchemski16}
are satisfied (w.r.t\ the Markov chain~$\Ceb$ and the limiting
process~$\CeLb$). Our statement thus follows from Theorems~1 and~2
there\footnote{Strictly speaking, the results are only stated when there
  is no killing, that is $\Cf(\Pa)=0$, but as mentioned by the authors
  \cite[p.~2562,~§2]{Kortchemski16}, they can be extended using the same
  techniques to the case where some killing is involved.}.
\end{prf}

Heading now toward pathwise and optional many-to-one formulae, we first
set up some notation. Let~$\Hs\subseteq\NN$ be a fixed subset of
states, and let~$\ell\in\partial\HT$ refer to an infinite word
over~$\NN$, which we see as a branch of~$\HT$. For every
$u\in\HT\cup\partial\HT$ and every $k\ge0$, set
\begin{equation*}
  \Cet_u(k)\defeq\Ce_{u[k]}(k-\Bt_{u[k]}),
\end{equation*}
where~$u[k]$ is the youngest ancestor~$v$ of~$u$ with $\Bt_v\le k$, and
write $\St[\Hs]_u\defeq\inf\{k\ge0\colon\Cet_u(k)\in\Hs\}$ for the
first hitting time of~$\Hs$ by~$\Cet_u$. Let also
$\Stb[\Hs]\defeq\inf\{k\ge0\colon\Ceb(k)\in\Hs\}$. Now, imagine that
once a particle hits~$\Hs$, it is stopped and thus has no further
progeny. The state when all particles have hit~$\Hs$ in finite time is
$\Sf[\Hs]_u\defeq\Cet_u(\St[\Hs]_u),\,u\in\HT_\Hs$, where
$\HT_\Hs\defeq
  \{u\in\HT\colon\ell[\St[\Hs]_\ell]=u
    \text{ for some }\ell\in\partial\HT\text{ with }
    \St[\Hs]_\ell<\infty\}$.

\begin{lem}[Many-to-one formula]{lem:Dirichlet}
\begin{enumerate}[label=(\roman*)]
  \item\label{itm:many-to-one} For every $n\in\NN$, every $k\ge0$,
    and every $f\colon\NN^{k+1}\to\RR_+$,
\begin{equation*}
 \Ex^{(n)}\!\left[\sum_{u\in\HT}
   \bigl(\Ce_u(k-\Bt_u)\bigr)^\Pa
   f\!\left(\Cet_u(i)\colon i\le k\right)
   \right]
 =\,n^\Pa\,
   \Ex^{(n)}\!\left[f\!\left(\Ceb(i)\colon i\le k\right)\!;\,
     \Ltb>k\right]\!.
\end{equation*}
  \item For every $n\in\NN$, every $\Hs\subseteq\NN$, and every
    $f\colon\ZZ_+\times\NN\to\RR_+$,
\begin{equation*}
 \Ex^{(n)}\!\left[\sum_{u\in\HT_\Hs}\!\bigl(\Sf[\Hs]_u\bigr)^\Pa
   f\bigl(\St[\Hs]_u,\Sf[\Hs]_u\bigr)\right]
 =\,n^\Pa\,
   \Ex^{(n)}\!\left[f\bigl(\Stb[\Hs],\Ceb\bigl(\Stb[\Hs]\bigl)\bigr);\,
     \Ltb>\Stb[\Hs]\right]\!.
\end{equation*}
\end{enumerate}
\end{lem}
\begin{prf}{prf:Dirichlet}
  (i)\enspace The proof is classical (see e.g.\ 
  \cite[Theorem~1.1]{Shi15}) and proceeds by induction on~$k$. The
  identity clearly holds for~$k=0$. Using~\eqref{eq:Weierstrass}
  together with the branching property at time~$k$,
\begin{flalign*}
\qquad&\Ex^{(n)}\!\left[\sum_{u\in\HT}
   \bigl(\Ce_u(k+1-\Bt_u)\bigr)^\Pa
   f\!\left(\Cet_u(i)\colon i\le k+1\right)
   \Big|\;\Cet_u(i)=x_{u,i},\,i\le k\right]&&\\[.4em]
&\qquad=\,\sum_{u\in\HT}
   \sum_{\vphantom{HT}m=0}^\infty\Tk_{x_{u,k},m}\left(m^\Pa
       f\bigl(x_{u,0},\ldots,x_{u,k},m\bigr)
     +\bigl(x_{u,k}-m\bigr)_+^\Pa
       f\bigl(x_{u,0},\ldots,x_{u,k},x_{u,k}-m\bigr)\right)\\[.4em]
&\qquad=\,\sum_{u\in\HT}\bigl(x_{u,k}\bigr)^\Pa\,
   \Ex^{(x_{u,k)}}\!\left[
    f\bigl(x_{u,0},\ldots,x_{u,k},\Ceb(1)\bigr);\,\Ceb(1)\neq0\right]\!.
\intertext{%
By taking expectations on both sides and applying the
induction hypothesis with the function
$\widetilde f(x_0,\ldots,x_k)
  \defeq\Ex^{(x_k)}[f(x_0,\ldots,x_k,\Ceb(1));\,\Ceb(1)\neq0]$
on the one hand, and by applying the Markov property of~$\Ceb$ at
time~$k$ on the other hand, we derive the identity at time~$k+1$:}
&\Ex^{(n)}\!\left[\sum_{u\in\HT}
  \bigl(\Ce_u(k+1-\Bt_u)\bigr)^\Pa
   f\!\left(\Cet_u(i)\colon i\le k+1\right)\right]
=\,n^\Pa\,
   \Ex^{(n)}\!\left[\widetilde f\!\left(\Ceb(i)\colon i\le k\right)\!;\,
     \Ltb>k\right]&\\[.4em]
&\phantom{\Ex^{(n)}\!\left[\sum_{u\in\HT}
  \bigl(\Ce_u(k+1-\Bt_u)\bigr)^\Pa
   f\!\left(\Cet_u(i)\colon i\le k+1\right)\right]}
=\,n^\Pa\,
    \Ex^{(n)}\!\left[f\!\left(\Ceb(i)\colon i\le k+1\right)\!;\,
     \Ltb>k+1\right]\!.&%
\end{flalign*}
(ii)\enspace For every $k\ge0$ and every $x_0,\ldots,x_k\in\NN$, we set
  $f^\Hs_k(x_0,\ldots,x_k)\defeq
    \II_{\{x_0\notin\Hs,\ldots,x_{k-1}\notin\Hs,x_k\in\Hs\}}f(k,x_k)$.
  Then
\begin{align*}
 \Ex^{(n)}\!\left[\sum_{u\in\HT_\Hs}\!\bigl(\Sf[\Hs]_u\bigr)^\Pa
   f\bigl(\St[\Hs]_u,\Sf[\Hs]_u\bigr)\right]
 &=\,\sum_{k=0}^\infty
   \Ex^{(n)}\!\left[\sum_{u\in\HT}\bigl(\Ce_u(k-b_u)\bigr)^\Pa
   f^\Hs_k\!\left(\Cet_u(i)\colon i\le k\right)
   \right]\\[.4em]
 &=\,n^\Pa\,\sum_{k=0}^\infty
   \Ex^{(n)}\!\left[f\bigl(k,\Ceb(k)\bigr);\,\Stb[\Hs]=k;\,
     \Ltb>k\right]\\[.4em]
 &=\,n^\Pa\,
   \Ex^{(n)}\!\left[f\bigl(\Stb[\Hs],\Ceb\bigl(\Stb[\Hs]\bigl)\bigr);\,
     \Ltb>\Stb[\Hs]\right]\!,
\end{align*}
by~\hyperref[itm:many-to-one]{(i)} and the monotone convergence theorem.
\end{prf}

We now combine \cref{pro:Descartes} and \cref{lem:Dirichlet} to
derive the following counterpart of \cite[Lemma~14]{Kortchemski18} that
we will apply in the next two sections. Consider the hitting set
$\Hs\defeq\{1,\ldots,\lfloor n\eps\rfloor\}$ and denote by
$\Sf_u\defeq\Sf[\Hs]_u,\,u\in\HT^{\le n\eps}\defeq\HT_\Hs,$ the
population of particles stopped below~$n\eps$.
\begin{cor}{cor:Ramanujan}
  We have
  \begin{equation*}
    \mathop{\vphantom{\limsup}\lim}_{\eps\to0}\,
    \limsup_{\vphantom{\eps}n\to\infty}\;
    n^{-\Pa}\,
    \Ex^{(n)}\!\left[\sum_{u\in\HT^{\le n\eps}}\!\!\!\!
      \bigl(\Sf_u\bigr)^\Pa\right]
    =\,0.
  \end{equation*}
\end{cor}
\begin{prf}{prf:Ramanujan}
  By \cref{lem:Dirichlet},
\begin{align*}
  n^{-\Pa}\,
    \Ex^{(n)}\!\left[\sum_{u\in\HT^{\le n\eps}}\!\!\!\!
      \bigl(\Sf_u\bigr)^\Pa\right]
    &=\,\Pr^{(n)}\!\left(\Ltb>\Stb\right)\!,
\intertext{where $\Stb\defeq\inf\{k\ge0\colon\Ceb(k)\le n\eps\}$. Thus,
  if~$\LtLb$ is the lifetime of~$\CeLb$ and
  $\StL\defeq\inf\{t\ge0\colon\CeLb(t)\le\eps\}$, then by
  \cref{pro:Descartes} and the continuous mapping theorem,}
    \limsup_{n\to\infty}\;
    n^{-\Pa}\,
    \Ex^{(n)}\!\left[\sum_{u\in\HT^{\le n\eps}}\!\!\!
      \bigl(\Sf_u\bigr)^\Pa\right]
    &\le\,\Pr\!\left(\LtLb>\StL\right)\!,
  \end{align*}
  which tends to~$0$ as $\eps\to0$ (because $\LtLb<\infty$
  and $\CeLb(\LtLb-)>0$, $\Pr$-almost surely).
\end{prf}

\section{Proof of \cref*{thm:Hilbert}}\label{sec:Hilbert}
We prove \cref{thm:Hilbert} by combining \cref{lem:Euler} with the next
two ``tightness'' properties. We suppose
that~\cref{assump:H1,assump:H2,assump:H3} hold and recall that
$\HTh\subset\HT$ refers to the set of words over $\{1,\ldots,h\}$
with length at most~$h$.
\begin{lem}{lem:Jacobi1}
  For every $\delta>0$,
\begin{equation*}
  \mathop{\vphantom{\limsup}\lim}_{h\to\infty}\;
    \Pr\!\left(\sup_{\vphantom{h}t\ge0}
      \sum_{u\in\HT\setminus\HTh}\!\!\!\!
      \bigl(\CeL_u(t-\BtL_u)\bigr)^\Pb\!>\delta\right)=\;0.
\end{equation*}
\end{lem}
\begin{prf}{lem:Jacobi1}
  This was already derived in \cite[Lemma~20]{Kortchemski18},
  and results from the following fact \cite[Corollary~4]{Bertoin17}:
  \begin{equation*}
    \Ex\!\left[\sum_{u\in\HT}\sup_{t\ge0}\;
      \bigl(\CeL_u(t-\BtL_u)\bigr)^q\right]<\,\infty
      \quad\text{for}\enspace\Cf(q)<0.
  \end{equation*}
\end{prf}
\begin{lem}{lem:Jacobi2}
  If~$\Ba$ is sufficiently large, then for every $\delta>0$,
\begin{equation*}
  \mathop{\vphantom{\limsup}\lim}_{h\to\infty}
  \limsup_{\vphantom{h}n\to\infty}\;
    \Pr^{(n)}\!\left(\sup_{\vphantom{h}k\ge0}
      \sum_{u\in\HT\setminus\HTh}\!\!\!\!
      \bigl(\Ce_u(k-\Bt_u)\bigr)^\Pb>\delta n^\Pb\right)
  =\;0.
\end{equation*}
\end{lem}
\begin{prf}{lem:Jacobi2}
  Let us first take~$\Pa<\Pb$ and~$\Ba$ as in \cref{sec:Fermat21}. As in
  the proof of~\cite[Lemma~22]{Kortchemski18} and by definition
  of~$\Good$ in~\cref{sec:Euler}, we claim that each particle in
  $\{\Ce_u(k-\Bt_u)\colon u\in\HT\setminus\Good\}$ has an ancestor
  with size at birth smaller than~$n\eps$. Thanks to the branching
  property, we may therefore consider that these particles derive from a
  system that has first been ``frozen'' below the level~$n\eps$, that
  is, with the notations of \cref{sec:Fermat21}, from a particle system
  having $\Sf_u,\,u\in\HT^{\le n\eps},$ as initial population. Hence, by
  \cref{lem:Abel} and Doob's maximal inequality,
  \begin{align*}
    \Pr^{(n)}\!\left(\sup_{k\ge0}
      \sum_{u\in\HT\setminus\Good}\!\!\!\!
        \bigl(\Ce_u(k-\Bt_u)\bigr)^\Pb>\delta n^\Pb\right)
      &\le\,\frac1{\delta^{\Pa/\Pb}\,n^\Pa}\,%
        \Ex^{(n)}\!\left[\sum_{u\in\HT^{\le n\eps}}\!\!\!
        \bigl(\Sf_u\bigr)^\Pa\right]
  \end{align*}
  (bounding from above the $\leb[\Pb]$-norm by the $\leb[\Pa]$-norm).
   We conclude by \cref{cor:Ramanujan} and \cref{lem:Neumann}.
\end{prf}
\begin{prf}[Proof of~\cref*{thm:Hilbert}]{thm:Hilbert}
  From \cref{lem:Euler,lem:Jacobi1,lem:Jacobi2}, we deduce the
  convergence in distribution
\begin{equation*}
  \left(\frac{\Ce_u^{(n)}(\lfloor\Ss_nt\rfloor-\Bt^{(n)}_u)}n\colon
     u\in\HT\right)_{t\ge0}
  \convlaw\,\bigl(Y_u(t-\BtL_u)\colon u\in\HT\bigr)_{t\ge0},
\end{equation*}
  in the space~$\Sk(\co{0}{\infty},\leb[\Pb](\HT))$ of
  $\leb[\Pb](\HT)$-valued càdlàg functions on~$\co{0}{\infty}$, where
\begin{equation*}
  \leb[\Pb](\HT)\defeq\left\{\bs x\defeq(x_u\colon
     u\in\HT)\colon\sum_{u\in\HT}(x_u)^\Pb<\infty\right\}\!.
\end{equation*}
  Since for $q\ge1$, rearranging sequences in the non-increasing order
  does not increase their $q$-distance \cite[Theorem~3.5]{Lieb01}, the
  convergence in~$\leb[\Pb](\HT)$ implies that
  in~$\leb[q\downarrow],\,q\ge1\vee\Pb$.
\end{prf}

\section{Proof of~\cref*{thm:Grothendieck}}\label{sec:Grothendieck}
Similarly to the previous section, by~\cref{lem:Galois} the proof
of~\cref{thm:Grothendieck} is complete once we have established that
\begin{gather}
  \lim_{h\to\infty}\Pr\Bigl(\dRT\bigl(\FTL,\FTL_h)>\delta\Bigr)
    \,=\,0,\notag
\shortintertext{and}
  \mathop{\vphantom{\limsup}\lim}_{h\to\infty}
    \limsup_{\vphantom{h}n\to\infty}\;
      \Pr\!\left(\dRT\Bigl(\FT^{(n)},
        \FT^{(n)}_h\Bigr)>\delta\Ss_n\right)
    =\,0,
\label{eq:Ramanujan}
\end{gather}
for all~$\delta>0$. The first display is clear since the tree~$\FTL$ is
compact. The second will be a consequence of the following counterpart of
\cite[Conjecture~1]{Kortchemski18}:
\begin{lem}{lem:Cantor}
  Suppose \labelcref{assump:H1,assump:H2,assump:H3}, and $\Pb>\Si$.
  Then for every $q<\Pb$, and for~$\Ba$ sufficiently large,
  \begin{equation*}
    \sup_{n\in\NN}\;\Ex\!\left[
      \left(\frac{\height\bigl(\FT^{(n)}\bigr)}{%
      \Ss_n}\right)^{\!q/\Si\,}\right]
      <\,\infty,
  \end{equation*}
  where
  $\height\bigl(\FT^{(n)}\bigr)\defeq\sup_{x\in\FT^{(n)}}\FTm(\es,x)$
  is the height of the tree~$\FT^{(n)}$.
\end{lem}
The proof of \cref{lem:Cantor} involves martingale arguments.
Prior to writing it, we need a preparatory lemma. Let us
define
\begin{equation*}
  \Cft_n(q)\defeq
    \Ss_n\sum_{m=1}^\infty\Tk_{n,m}
      \left[\left(\frac{\Ss_m}{\Ss_n}\right)^{\!q/\Si}-1
        +\left(\frac{\Ss_{n-m}}{\Ss_n}\right)^{\!q/\Si\,}
       \right]\!,
\end{equation*}
which slightly differs from~$\Cf_n(q)$ to the extent that we have
replaced the map $m\mapsto m^q$ by the $q$-regularly-varying sequence
$\Ssr{q}(m)\defeq\Ss_m^{q/\Si},\,m\in\NN$ (for convenience, we have set
$\Ss_m\defeq0,\,m\le0$). Of~course, $\Cft_n=\Cf_n$ if
$\Ss_m=m^\Si$ for every $m\in\NN$.
\begin{lem}{lem:Cauchy}
  Suppose $\Pb>\Si$. Then we can find $\Pa\in\oo{0}{\Pb}$ such that, for every
  $q\in\co{\Pa}{\Pb}$,
  \begin{equation*}
    \lim_{n\to\infty}\Cft_n(q)=\Cf(q)<0.
  \end{equation*}
\end{lem}
\begin{prf}{prf:Cauchy}
  We will more generally show that for every $q$-regularly-varying
  sequence~$(r_n)$,
  \begin{equation*}
    \left|\Ss_n\sum_{m=1}^\infty\Tk_{n,m}
      \left[\frac{r_m}{r_n}-\left(\frac mn\right)^{\!q\,}\right]\right|
   +\left|\Ss_n\sum_{m=1}^{n-1}\Tk_{n,m}\left[\frac{r_{n-m}}{r_n}
          -\left(1-\frac mn\right)^{\!q\,}\right]\right|
          \conv\,0,
  \end{equation*}
  provided $q<\Pb$ is close enough to~$\Pb$.
  Denoting by~$\Jm_n$ the law of $\log(\Ce^{(n)}(1)/n)$, we observe that
  \begin{equation*}
    \Ss_n\sum_{m=1}^\infty\Tk_{n,m}
      \left[\frac{r_m}{r_n}-\left(\frac mn\right)^{\!q\,}\right]
    =\Ss_n\int_{-\infty}^\infty
      \left[\left(\frac{r_{ne^x}}{r_n}\right)
        -e^{qx}\right]\Jm_n(\dd x),
  \end{equation*}
  which, by repeating the arguments in
  \cite[Proof of Lemma~4.9]{Kortchemski16}, tends to~$0$ as
  $n\to\infty$. Next, an appeal to Potter's bounds
  \cite[Theorem~1.5.6]{Bingham87}
  shows that for every $c>1$ and $\delta>0$ arbitrary small,
  \begin{equation*}
      \frac1c\left(\frac mn\right)^{\!q+\delta}
      \le\frac{r_m}{r_n}
      \le\,c\left(\frac mn\right)^{\!q-\delta}
  \end{equation*}
  whenever $m<n$ are sufficiently large. Thus, recalling that
  $\Le_n(q)\to\Le(q)$ and $\Cf_n(q)\to\Cf(q)$ for every $q$ in some
  left-neighbourhood of~$\Pb$, we have
  \begin{align*}
    \liminf_{n\to\infty}\;\Ss_n\sum_{m=1}^{n-1}\Tk_{n,m}\left[
      \frac{r_{n-m}}{r_n}-\left(1-\frac mn\right)^{\!q\,}\right]
    &\ge\,\frac1c\,\Bigl(\Cf(q+\delta)-\Le(q+\delta)\Bigr)
      -\Bigl(\Cf(q)-\Le(q)\Bigr),\\
    \intertext{and}
    \limsup_{n\to\infty}\;\Ss_n\sum_{m=1}^{n-1}\Tk_{n,m}\left[
      \frac{r_{n-m}}{r_n}-\left(1-\frac mn\right)^{\!q\,}\right]
    &\le\,\frac1c\,\Bigl(\Cf(q-\delta)-\Le(q-\delta)\Bigr)
      -\Bigl(\Cf(q)-\Le(q)\Bigr).
  \end{align*}
  We conclude by letting $c\to1$ and $\delta\to0$.
\end{prf}

\noindent
We can now prove \cref*{lem:Cantor}.
\begin{prf}[Proof of \cref*{lem:Cantor}]{lem:Cantor}
  We shall rely on a Foster-type technique close to the machinery
  developed in~\cite{Aspandiiarov98}; see in particular the proof of
  Theorem~2' there. First, observe that~$\height(\FT)$ is distributed
  like the extinction time~$\Et$ of~$\GF$:
  \begin{equation*}
    \height\bigl(\FT^{(n)}\bigr)
    \eqlaw\,\sup_{u\in\HT}\;\Et^{(n)}_u\,
    \eqdef\,\Et^{(n)}.
  \end{equation*}
  Fix $q\in\oo{\Si}{\Pb}$ arbitrary close to~$\Pb$ and set
  $r\defeq q/\Si$. By \cref{lem:Cauchy}, suppose~$\Ba$ large enough so
  that $\Cft_m(q)<0$ for every $m>\Ba$. It is easy to see as in the
  proof of \cref{lem:Abel} that the process
\begin{equation*}
\Gamma(k)\defeq\sum_{u\in\HT}
   \Ssr{q}\bigl(\Ce_u(k-\Bt_u)\bigr),\qquad k\ge0,
\end{equation*}
is a supermartingale under~$\Pr^{(n)}$ (with respect to the natural
filtration $(\mathcal F_k)_{k\ge0}$ of~$\GF$): indeed, for
$\GF(k)=(x_i\colon i\in I)$,
\begin{equation*}
  \Ex^{(n)}\!\left[\Gamma(k+1)-\Gamma(k)\;\Big|\;
    \mathcal F_k\right]
    =\,\sum_{i\in I}\Cft_{x_i}(q)\,\Ssr{q-\Si}(x_i),
\end{equation*}
where the right-hand side is (strictly) negative on the event
$\{\Et>k\}=\{\exists i\in I\colon x_i>\Ba\}$. We will more precisely
show the existence of $\eta>0$ sufficiently small such that the process
\begin{equation*}
  G(k)\defeq\left(\Gamma(k)^{1/r}+\eta\,\bigl(\Et\wedge k
    \bigr)\right)^{\!r\!},\qquad k\ge0,
\end{equation*}
is a $(\mathcal F_k)_{k\ge0}$-supermartingale under $\Pr^{(n)}$, for any
$n\in\NN$. Then, the result will be readily obtained from
$\eta^r\,\Ex^{(n)}[(\Et\wedge k)^r]\le\Ex^{(n)}[G(k)]\le\Ex^{(n)}[G(0)]
=\Ssr{q}(n)=\Ss_n^r$ and an appeal to Fatou's lemma.

On the one hand, we have
\begin{equation*}
  \sigma\defeq\sum_{i\in I}\Ssr{q-\Si}(x_i)
  \ge\left(\sum_{i\in I}\Ssr{q}(x_i)\right)^{\!1-\Si/q}
\end{equation*}
because
\begin{equation*}
  \frac{\Ssr{q}(x_i)}{\sigma^{q/(q-\Si)}}\,=
  \left(\frac{\Ssr{q-\Si}(x_i)}\sigma\right)^{\!q/(q-\Si)}
  \le\,\frac{\Ssr{q-\Si}(x_i)}{\sigma},
\end{equation*}
where $q/(q-\Si)>1$ and the right-hand side sums to~$1$ as~$i$ ranges
over~$I$. Then, if we let $\eta>0$ sufficiently small such that
$\Cft_m(q)\le -r\eta$ for every $m>\Ba$, we deduce that
\begin{equation*}
  \Ex^{(n)}\!\left[\Gamma(k+1)\;\Big|\;\mathcal F_k\right]
  \le\,\Gamma(k)\left(1-r\eta\,\Gamma(k)^{-\Si/q}
    \,\II_{\{\Et>k\}}\right)\!.
\end{equation*}
Raising this to the power~$1/r=\Si/q$ yields
\begin{equation}
  \Ex^{(n)}\!\left[\Gamma(k+1)\;\Big|\;\mathcal F_k\right]^{1/r}
  \le\,\Gamma(k)^{1/r}\left(1-\eta\,\Gamma(k)^{-\Si/q}
    \,\II_{\{\Et>k\}}\right)
  =\,\Gamma(k)^{1/r}-\eta\,\II_{\{\Et>k\}},
  \label{eq:Weyl}
\end{equation}
by concavity of $x\mapsto x^{1/r}$. On the other hand, the
supermartingale property also implies that $(\Gamma(k+1)^{1/r}+a)^r$ is
integrable for every constant $a>0$; we may thus apply the generalized
triangle inequality \cite[Lemma~1]{Aspandiiarov98} with the positive,
convex increasing function $x\mapsto x^r$, the positive random variable
$\Gamma(k+1)^{1/r}$, and the probability
$\Pr^{(n)}(\;\cdot\mid\mathcal F_k)$ (under which $\Et\wedge(k+1)$ can
be seen as a positive constant):
\begin{equation*}
  \Ex^{(n)}\!\left[\left(\Gamma(k+1)^{1/r}+\eta\,\bigl(\Et\wedge(k+1)
    \bigr)\right)^r\;\Big|\;\mathcal F_k\right]^{1/r}
  \le\;\Ex^{(n)}\bigl[\Gamma(k+1)\;\big|\;\mathcal F_k\bigr]^{1/r}
    +\eta\,\bigl(\Et\wedge(k+1)\bigr).
\end{equation*}
Reporting~\eqref{eq:Weyl} shows as desired that $(G(k)\colon k\ge0)$
is a supermartingale.
\end{prf}

\noindent
We are finally ready to derive~\eqref{eq:Ramanujan} and complete the proof
of \cref{thm:Grothendieck}.
\begin{prf}[Proof of~\labelcref*{eq:Ramanujan}]{eq:Ramanujan}
We start as in the proof of \cref{thm:Hilbert}: thanks
to~\cref{lem:Neumann} and the branching property, with high probability
as $h\to\infty$, the connected components of
$\FT^{(n)}\setminus\FT^{(n)}_h$ are included in independent copies
of~$\FT$ stemming from the population $\Sf_u,\,u\in\HT^{\le n\eps},$ of
particles frozen below~$n\eps$. Specifically,
\begin{equation*}
  \Pr^{(n)}\!\left(\dRT\Bigl(\FT^{(n)},
    \FT^{(n)}_h\Bigr)>\delta\Ss_n\right)
  \le\,\Pr^{(n)}\!\left(\Good\not\subseteq\HTh\right)
    +\Ex^{(n)}\!\left[
      \sum_{u\in\HT^{\le n\eps}}\!\!\!
      \Pr^{(\Sf_u)}\Bigl(\height(\FT)>\delta\Ss_n\Bigr)\right]\!.
\end{equation*}
Now, take $\Pa<q<\Pb$ and~$\Ba$ large enough so that both
\cref{lem:Cantor} and the results of \cref{sec:Fermat21} hold. So, there
exists a constant $C>0$ such that
\begin{align*}
\Ex^{(n)}\!\left[
\sum_{u\in\HT^{\le n\eps}}\!\!\!
  \Pr^{(\Sf_u)}\Bigl(\height(\FT)>\delta\Ss_n\Bigr)\right]
  &=\,\Ex^{(n)}\!\left[\sum_{u\in\HT^{\le n\eps}}\!\!\!
    \Pr^{(\Sf_u)}\!\left(\height(\FT)>\delta\Ss_{\Sf_u}
      \frac{\Ss_n}{\Ss_{\Sf_u}}\right)\right]\\[.4em]
  &\le\,C\,\Ex^{(n)}\!\left[\sum_{u\in\HT^{\le n\eps}}\!\!
    \left(\frac{\Ss_{\Sf_u}}{\Ss_n}\right)^{\!q/\Si\,}\right]\!.
\end{align*}
But we know, thanks to another application of Potter's bounds, that we
may find $c>0$ such that
\begin{equation*}
\left(\frac{\Ss_m}{\Ss_n}
    \right)^{q/\Si}\le c\left(\frac mn\right)^{\!\Pa\,}\!,
\end{equation*}
whenever~$n$ is sufficiently large and $m\le n$. Since $\Sf_u\le n$ (for
$0<\eps<1$), we can again conclude by \cref{cor:Ramanujan} and
\cref{lem:Neumann}.
\end{prf}

\section*{Acknowledgments}
The author would like to thank Jean Bertoin for his constant suggestions
and guidance, and also Bastien Mallein for helpful discussions and
comments. Thanks are also directed to an anonymous referee for their careful
reading.

\end{document}